\documentclass[11pt]{article}
\usepackage{amsfonts}
\usepackage{latexsym}
\usepackage{amsmath}
\usepackage{amssymb}
\usepackage{amsthm}
\makeatletter \oddsidemargin -.1in \evensidemargin -.1in
\newcommand{\singlespacing}{\let\CS=\@currsize\renewcommand{\baselinestretch}{1}\tiny \CS}
\newcommand{\doublespacing}{\let\CS=\@currsize\renewcommand{\baselinestretch}{1.35}\tiny \CS}
\date{}
\newcommand{\be}{\begin{equation}}
\newcommand{\ee}{\end{equation}}
\newcommand{\bea}{\begin{eqnarray}}
\newcommand{\eea}{\end{eqnarray}}
\newcommand{\nn}{\nonumber}
\newcommand{\bee}{\begin{eqnarray*}}
\newcommand{\eee}{\end{eqnarray*}}

\newtheorem{e1}{Example}[section]
\newtheorem{t1}{Theorem}[section]
\newtheorem{d1}{Definition}[section]

\newtheorem{l1}{Lemma}[section]
\newtheorem{p1}{Proposition}[section]

\begin{document}
\title{Controlled $g$-frames in Hilbert $C^*$-modules}
\author{N. K. Sahu\thanks{E-mail:~nabindaiict@gmail.com}
\\
\small \it{{*}Dhirubhai Ambani Institute of Information and
Communication Technology, Gandhinagar, India}}

\date{}
\maketitle \setlength{\baselineskip}{18pt} \noindent 
\begin{abstract}
To improve the numerical efficiency of iterative algorithms for
inverting the frame operator, the controlled frame was introduced
by Balazs et al. \cite{Balazs}, and has since been given more
importance. In this paper, we introduce the concept of controlled
g-frames in Hilbert $C^{*}$-modules. We establish the equivalent
condition for controlled $g$-frame using operator theoretic
approach. We investigate some operator theoretic characterizations
of controlled $g$-frames and controlled $g$-Bessel sequences. We
also established the relationship between $g$-frames and
controlled $g$-frames in Hilbert $C^{*}$-modules. At the end we
prove some perturbation results on controlled $g$-frames.
\end{abstract}



\section{Introduction}
It was Duffin and Schaeffer \cite{Duffin} in 1952 who initiated
the notion of frames while studying nonharmonic Fourier series.
After a long gap, in 1986, Daubechies et al. \cite{Daubechies}
reintroduced the same notion and developed the theory of frames.
In general, frame is nothing but a spanning set and what makes it
interesting is the redundance of additional vectors than those in
the basis. Due to its redundancy it becomes more applicable not
only in theoretical point of view but also in various kinds of
applications. Due to their rich structure the subject draws the
attention of many mathematicians,
 physicists and engineers since it is largely applicable in signal processing \cite{Ferreira}, image processing \cite{Cand}, coding and communications \cite{Strohmer},
 sampling \cite{Eldar, Eldar Y.}, numerical analysis, filter theory \cite{olcskei}. Now a days it is used in compressive sensing, data analysis and other areas.

Hilbert $C^{*}$-module is a wide category between Hilbert space
and Banach space.  It was Frank and Larson \cite{Frank}, who
initiated the theory of frames in Hilbert $C^{*}$-modules. For
more details of frames in Hilbert $C^{*}$-modules one may refer to
the Doctoral Dissertation \cite{Jing}, Han et al. \cite{Han1}, Han
et al. \cite{Han2} and the references there in. The notion of
g-frame or generalized frame in Hilbert $C^{*}$-module is
introduced by Sun \cite{Sun}. For more on $g$-frames one can refer
to Khosravi and Khosravi \cite{Khosravi}, Fu and Zhang \cite{Fu},
Li and Leng \cite{Li}. To improve the numerical efficiency of
iterative algorithms for inverting the frame operator, controlled
frame was introduced by Balazs et al. \cite{Balazs} in Hilbert
spaces. Recently, Kouchi
 and Rahimi \cite{Rashidi} introduced Controlled frames in Hilbert $C^{*}$-modules.
  Motivated from the above literature, we introduce the notion of controlled g-frame in Hilbert $C^{*}$-modules.\\
\section{Preliminaries}
\noindent Let us briefly recall some definitions and basic
properties of Hilbert $C^{*}$ modules. Hilbert $C^{*}$-modules are
generalization of Hilbert spaces by allowing the inner product to
take values in a $C^{*}$-algebra rather than the usual fields
 $\mathbb{R}$ or $\mathbb{C}$.\\
Let $\mathcal{A}$ be a unital $C^{*}$-algebra. The Hilbert
$C^{*}$-module or Hilbert $\mathcal{A}$-module is defined as
follows:


\begin{d1}\rm   Let $\mathcal{H}$ be a left $\mathcal{A}$-module such that the linear structure of $\mathcal{A}$ and $\mathcal{H}$ are compatible;
 $\mathcal{H}$ is called a pre-Hilbert $\mathcal{A}$-module if $\mathcal{H}$ is equipped with an $\mathcal{A}$-valued inner product
  $\langle \cdot, \cdot \rangle: \mathcal{H} \times \mathcal{H} \to \mathcal{A}$ such that \\
(i) $\langle f, f \rangle \geq 0$ $, \forall  f \in \mathcal{H}$, and $\langle f, f \rangle = 0$ if and only if $f = 0$\\
(ii) $\langle f, g \rangle  = \langle g, f \rangle^* $\\
(iii) $\langle af+g,h \rangle  = a\langle f, h \rangle + \langle g, h \rangle $ for all $f,g,h \in \mathcal{H}$ and $a \in \mathcal{A}$.
\end{d1}
\noindent For every $f \in \mathcal{H}$, the norm is defined as $\|f\| = \| \langle f, f \rangle \|^{\frac{1}{2}}$.\\
If $\mathcal{H}$ is complete with respect to the norm, it is called a Hilbert $\mathcal{A}$-module or a Hilbert $C^{*}$-module over  $\mathcal{A}$.\\
Initially, it was Frank and Larson \cite{Frank} who introduced the
notion of frames in Hilbert $C^{*}$-modules. Their definition is
the following:

\begin{d1}\rm  \cite{Frank} A set of elements $\{f_j\}_{j \in \mathbb{J}}$ in a Hilbert $C^{*}$-module $\mathcal{H}$ over a unital $C^{*}$-algebra $\mathcal{A}$,
is said to be a frame if there exist two constants $C, D > 0$ such
that \bee
        C \langle f,f \rangle \leq \sum_{j \in \mathbb{J}} \langle f,f_j \rangle \langle f_j,f \rangle \leq D\langle f,f \rangle ,
        \forall f \in \mathcal{H}.
\eee
\end{d1}
After the introduction of frames in Hilbert $C^{*}$-modules, a lot
of work on frame theory has been developed in Hilbert
$C^{*}$-modules.  The concept of g-frames in Hilbert
$C^{*}$-modules was introduced by Sun \cite{Sun}. The $g$-frame in
Hilbert $C^{*}$-modules is defined as follows:

\begin{d1}\label{d1}\rm \cite{Sun} Let $\mathcal{H}$ be a Hilbert $C^{*}$-module over $\mathcal{A}$, $\{\mathcal{H}_j\}_{j \in \mathbb{J}}$ be a
sequence of closed subspaces
 of $\mathcal{H}$.
A sequence $\{ \wedge_j \in
End_{A}^{*}(\mathcal{H},\mathcal{H}_j)\}$ is called a g-frame in
$\mathcal{H}$ with respect to $\{\mathcal{H}_j\}_{j \in
\mathbb{J}}$ if there exist positive constants $C, D$ such that
\bea \label{eqn17}
    C \langle f,f \rangle \leq \sum_{j \in \mathbb{J}} \langle \wedge_jf,\wedge_jf \rangle \leq D\langle f,f \rangle ,
\forall f \in \mathcal{H}. \eea
\end{d1}
The g-frame operator $S_{\Lambda}:H\rightarrow H$ is defined as
\bea \label{eqn18} S_{\Lambda}f=\sum_{j\in
J}\Lambda_{j}^{*}\Lambda_{j}f. \eea
To improve the numerical
efficiency of iterative algorithms to find the inverse of frame
operator, a new notion of frame was introduced, that is called as
 controlled frames. Controlled frames in Hilbert space is
 introduced by Balazs et al. \cite{Balazs} in 2010. Very recently,
 controlled frames in Hilbert $C^{*}$-modules is introduced by Kouchi
 and Rahimi \cite{Rashidi}. The controlled frame in Hilbert
 $C^{*}$-modules is defined as follow:
\begin{d1}\rm  \cite{Rashidi} Let $\mathcal{H}$ be a Hilbert $C^{*}$-module and $C \in GL(\mathcal{H})$. A family of vectors
 $\{f_j \in \mathcal{H}\}_{j \in \mathbb{J}}$
 is said to be a controlled frame in $\mathcal{H}$ or $C$-Controlled frame in $\mathcal{H}$ if there exist constants $m, M > 0$ such that
\bee
m \langle f,f \rangle \leq \sum_{j \in \mathbb{J}} \langle f,f_j \rangle \langle Cf_j,f \rangle \leq M\langle f,f \rangle , \forall f \in \mathcal{H}.
\eee
\end{d1}

In the next section we introduce the notion of Controlled g-frames
in Hilbert $C^{*}$-modules. We study several characterizations of
a controlled g-frame, equivalent formulation, its operator
theoretic behavior, its relationship with the frames etc.. In the
end we present some stability results on controlled $g$-frames.

\section{Controlled g-frame}
 Let $\mathcal{H}$ be a $C^{*}$-module over a unital $C^{*}$-algebra $\mathcal{A}$ with
 $\mathcal{A}$-valued inner product $\langle.,.\rangle$ and norm $\|.\|$.
 Let $\{\mathcal{H}_{j}\}_{j\in J}$ be a sequence of closed submodules of
 $\mathcal{H}$, where $J$ is any index set. Also let $GL^{+}(\mathcal{H})$ be the set
 of all positive bounded linear invertible operators on $\mathcal{H}$ with
 bounded inverse.
 \begin{d1}\rm
 Let $C,C^{\prime}\in GL^{+}(\mathcal{H})$. A sequence $\{\Lambda_{j}\in End_{\mathcal{A}}^{*}(\mathcal{H},\mathcal{H}_{j}):j\in J\}$ is
 said to be a $(C,C^{\prime})$-controlled g-frame for $\mathcal{H}$ with respect
 to $\{\mathcal{H}_{j}\}_{j\in J}$ if there exist constants $0<A\leq
 B<\infty$ such that
 \bea\label{eqn1}
 A\langle f,f\rangle\leq \sum_{j\in J}\langle \Lambda_{j}C f, \Lambda_{j}C^{\prime}
 f\rangle\leq B\langle f,f\rangle,~\forall f\in \mathcal{H}.
 \eea
 When $A=B$, the sequence $\{\Lambda_{j}\in End_{\mathcal{A}}^{*}(\mathcal{H},\mathcal{H}_{j}):j\in J\}$ is called $(C,C^{\prime})$-controlled tight
 g-frame, and when $A=B=1$, it is called a
 $(C,C^{\prime})$-controlled Parseval g-frame.
 \end{d1}

 \begin{d1}\rm
A sequence $\{\Lambda_{j}\in
End_{\mathcal{A}}^{*}(\mathcal{H},\mathcal{H}_{j}):j\in J\}$ is
 said to be a $(C,C^{\prime})$-controlled g-Bessel sequence for $\mathcal{H}$ with respect
 to $\{\mathcal{H}_{j}\}_{j\in J}$ if there exists constant $0<
 B<\infty$ such that
 \bea\label{eqn1.1}
 \sum_{j\in J}\langle \Lambda_{j}C f, \Lambda_{j}C^{\prime}
 f\rangle\leq B\langle f,f\rangle,~\forall f\in \mathcal{H}.
 \eea
\end{d1}

\begin{e1}\rm
Let $\mathcal{H}$ be an ordinary inner product space,
$J=\mathbb{N}$, and $\{e_{j}\}_{j=1}^{\infty}$ be an orthonormal
basis for Hilbert $\mathcal{C}$-module $\mathcal{H}$. We construct
$\mathcal{H}_{j}$ as
$\mathcal{H}_{j}=\overline{span}\{e_{1},e_{2},...,e_{j}\}$ for
each $j\in \mathbb{N}$.\\ Define
$\Lambda_{j}:\mathcal{H}\rightarrow \mathcal{H}_{j}$ by \bee
\Lambda_{j}f=\sum_{k=1}^{j}\big\langle
f,\frac{e_{j}}{\sqrt{j}}\big\rangle e_{k}.\eee The adjoint
operator $\Lambda_{j}^{*}:\mathcal{H}_{j}\rightarrow \mathcal{H}$
can be easily found as \bee
\Lambda_{j}^{*}(g)=\sum_{k=1}^{j}\big\langle
f,\frac{e_{k}}{\sqrt{j}}\big\rangle e_{j}.\eee Let us define
$Cf=2f$ and $C'f=\frac{1}{2}f$. Then for any $f\in \mathcal{H}$,
we can estimate \bee \Big\langle
\Lambda_{j}Cf,\Lambda_{j}C'f\Big\rangle&=&\Big\langle\sum_{k=1}^{j}\Big\langle
2f,\frac{e_{j}}{\sqrt{j}}\Big\rangle
e_{k},\sum_{k=1}^{j}\Big\langle
\frac{1}{2}f,\frac{e_{j}}{\sqrt{j}}\Big\rangle
e_{k}\Big\rangle\\
&=& \Big\langle 2f,\frac{e_{j}}{\sqrt{j}}\Big\rangle \Big\langle
\frac{1}{2}f,\frac{e_{j}}{\sqrt{j}}\Big\rangle^{*}\sum_{k=1}^{j}\|e_{k}\|^{2}\\
&=& \Big\langle 2f,\frac{e_{j}}{\sqrt{j}}\Big\rangle \Big\langle
\frac{1}{2}f,\frac{e_{j}}{\sqrt{j}}\Big\rangle^{*}~j\\
&=& \langle f,e_{j}\rangle\langle f,e_{j}\rangle^{*}.\eee
Therefore, for any $f\in \mathcal{H}$, \bee \sum_{j=1}^{\infty}
\Big\langle \Lambda_{j}Cf,\Lambda_{j}C'f\Big\rangle&=&
\sum_{j=1}^{\infty}\langle f,e_{j}\rangle\langle
f,e_{j}\rangle^{*}=\langle f,f\rangle.\eee This shows that
$\{\Lambda_{j}:j\in \mathbb{N}\}$ is a $(C,C')$-controlled
Parseval $g$-frame for $\mathcal{H}$ with respect to
$\{H_{j}\}_{j\in \mathbb{N}}$.
\end{e1}

Suppose that  $\{\Lambda_{j}\in
End_{\mathcal{A}}^{*}(\mathcal{H},\mathcal{H}_{j}):j\in J\}$ be a
$(C,C^{\prime})$-controlled g-frame for the Hilbert $C^{*}$-module
$\mathcal{H}$ with respect to $\{\mathcal{H}_{j}\}_{j\in J}$. The
bounded linear operator
$T_{(C,C^{\prime})}:l^{2}(\{\mathcal{H}_{j}\}_{j\in J})\rightarrow
\mathcal{H}$ defined by \bea\label{eqn2}
T_{(C,C^{\prime})}(\{g_{j}\}_{j\in J})=\sum_{j\in J}
(CC')^{\frac{1}{2}}\Lambda_{j}^{*}g_{j},~~\forall \{g_{j}\}_{j\in
J}\in l^{2}(\{\mathcal{H}_{j}\}_{j\in J})\eea
 is called the synthesis
operator for the $(C,C^{\prime})$-controlled g-frame
$\{\Lambda_{j}:j\in J\}$.\\ The adjoint operator
$T^{*}_{(C,C^{\prime})}: \mathcal{H}\rightarrow
l^{2}(\{\mathcal{H}_{j}\}_{j\in J})$ given by \bea\label{eqn3}
T^{*}_{(C,C^{\prime})}(f)=\big\{\Lambda_{j}(C^{\prime}C)^{\frac{1}{2}}f\big\}_{j\in
J} \eea is called the analysis operator for the
$(C,C^{\prime})$-controlled g-frame $\{\Lambda_{j}:j\in J\}$.\\
When $C$ and $C'$ commute with each other, and commute with the
operator $\Lambda_{j}^{*}\Lambda_{j}$ for each $j$, then the
$(C,C^{\prime})$-controlled g-frame operator
$S_{(C,C^{\prime})}:\mathcal{H}\rightarrow \mathcal{H}$ is defined
as \bea\label{eqn4} S_{(C,C^{\prime})}f &=&
T_{(C,C^{\prime})}T^{*}_{(C,C^{\prime})}f =\sum_{j\in J}
C^{\prime}\Lambda_{j}^{*}\Lambda_{j}C f.\eea For the above result
one is referred to Hua and Huang \cite{Hua}. So from now on we
assume that $C$ and $C'$ commute with each other, and commute with
the operator $\Lambda_{j}^{*}\Lambda_{j}$ for each $j$.

\begin{p1}\rm \label{pro1}
Let $\{\Lambda_{j}:j\in J\}$ be a $(C,C^{\prime})$-controlled
g-frame for the Hilbert $C^{*}$-module $\mathcal{H}$ with respect
to $\{\mathcal{H}_{j}\}_{j\in J}$. Then the
$(C,C^{\prime})$-controlled g-frame operator $S_{(C,C^{\prime})}$
is positive, self adjoint and invertible.
\end{p1}

\begin{proof}
The frame operator $S_{(C,C^{\prime})}$ for the
$(C,C^{\prime})$-controlled g-frame is
$S_{(C,C^{\prime})}f=\sum_{j\in J}
C^{\prime}\Lambda_{j}^{*}\Lambda_{j}Cf$. As $\{\Lambda_{j}:j\in
J\}$ is a $(C,C^{\prime})$-controlled g-frame, and from the
following identity, \bee
  \sum_{j\in J}\langle \Lambda_{j}C f, \Lambda_{j}C^{\prime}
 f\rangle=\big\langle\sum_{j\in J} C'\Lambda^{*}_{j}\Lambda_{j}C f,
 f\big\rangle=
 \big\langle S_{(C,C^{\prime})}f,f\big\rangle,
 \eee
 we clearly see that $S_{(C,C^{\prime})}$ is a positive operator.
 Also it is clearly bounded and linear. Again
\bee \big\langle S_{(C,C^{\prime})}f,g\big\rangle &=&
\big\langle\sum_{j\in J} C'\Lambda^{*}_{j}\Lambda_{j}C f,
 g\big\rangle\\
 &=&\sum_{j\in J}\big\langle C'\Lambda^{*}_{j}\Lambda_{j}C f,
 g\big\rangle\\
 &=&\sum_{j\in J}\big\langle f, C\Lambda^{*}_{j}\Lambda_{j}C'g
 \big\rangle=\sum_{j\in J}\big\langle f, S_{(C',C)}g
 \big\rangle.
 \eee
 Hence $S_{(C,C')}^{*}=S_{(C',C)}$. Also as $C$ and $C'$ commute
 with each other and commute with $\Lambda_{j}^{*}\Lambda_{j}$, we have
 $S_{(C,C')}=S_{(C',C)}$. So the controlled g-frame operator is
 self adjoint. Alternatively, this can also be directly obtained as $S_{(C,C')}$
 is a positive operator, and every positive operator is self
 adjoint.\\
 From the controlled g-frame identity we have
 \bee
 &&A\langle f,f\rangle\leq \big\langle S_{(C,C')}f,f\big\rangle\leq
 B\langle f,f\rangle\\
&\Rightarrow& A~Id_{\mathcal{H}}\leq S_{(C,C')}\leq
B~Id_{\mathcal{H}},
 \eee
 where $Id_{\mathcal{H}}$ is the identity operator in $\mathcal{H}$. Thus the
 controlled g-frame operator $S_{(C,C')}$ is invertible.
\end{proof}

\begin{l1} \label{l3.1} \rm \cite{Arambaic}  Let $\mathcal{A}$ be a $C^*$-algebra. Let $U$ and $V$ be two Hilbert $\mathcal{A}$-modules and $T\in End_{\mathcal{A}}^{*}(U,V)$.
Then the following statements are equivalent:
\begin{enumerate}
    \item $T$ is surjective.
    \item $T^*$ is bounded below with respect to norm i.e there exists $m > 0$ such that $\|T^{*}f\| \geq m \|f\|$ for all $f \in U$.
    \item  $T^*$ is bounded below with respect to inner product i.e there exists $m > 0$ such that $\langle T^{*}f,T^{*}f \rangle \geq m \langle f,f \rangle$ for all $f \in U$.
\end{enumerate}
\end{l1}

With the help of the above Lemma \ref{l3.1},  we establish an
equivalent definition of $(C,C')$-controlled g-frame.

\begin{t1}\label{t1}\rm
Let $\{\Lambda_{j}:j\in J\}\subset
End^{*}_{A}(\mathcal{H},\mathcal{H}_{j})$ and $\sum_{j\in
J}\langle \Lambda_{j}Cf,\Lambda_{j}C'f\rangle$ converge in norm
for any $f\in \mathcal{H}$. Then $\{\Lambda_{j}:j\in J\}$ is a
$(C,C')-$controlled g-frame for $\mathcal{H}$ with respect to
$\{\mathcal{H}_{j}\}_{j\in J}$ if and only if there exists
constants $A,B>0$ such that \bea A\|f\|^{2}\leq \|\sum_{j\in
J}\langle \Lambda_{j}Cf,\Lambda_{j}C'f\rangle\|\leq
B\|f\|^{2},~~\forall f\in \mathcal{H}.\eea
\end{t1}
\begin{proof}
Let $\{\Lambda_{j}:j\in J\}$ be a $(C,C')-$controlled g-frame for
$\mathcal{H}$ with respect to $\{\mathcal{H}_{j}\}_{j\in J}$ with
bound $A$ and $B$. Hence we have \bea\label{eqn5} A\langle
f,f\rangle\leq \sum_{j\in J}\langle \Lambda_{j}C f,
\Lambda_{j}C^{\prime}
 f\rangle\leq B\langle f,f\rangle,~\forall f\in \mathcal{H}.\eea
 Since  $\langle f,f\rangle\geq 0,~~\forall f\in \mathcal{H}$, then
 we can take the norm on the left, middle and right terms of
 the above inequality (\ref{eqn5}). Thus we have
 \bee &&\|A\langle f,f\rangle\|\leq \|\sum_{j\in
J}\langle \Lambda_{j}C f, \Lambda_{j}C^{\prime}
 f\rangle\|\leq \|B\langle f,f\rangle\|,~\forall f\in \mathcal{H}\\
 &\Rightarrow& A\|f\|^{2}\leq \|\sum_{j\in
J}\langle \Lambda_{j}C f, \Lambda_{j}C^{\prime}
 f\rangle\|\leq B\|f\|^{2},~~\forall f\in \mathcal{H}.
 \eee
 Conversely, suppose that
 \bea\label{eqn6} A\|f\|^{2}\leq \|\sum_{j\in
J}\langle \Lambda_{j}C f, \Lambda_{j}C^{\prime}
 f\rangle\|\leq B\|f\|^{2},~~\forall f\in \mathcal{H}.
 \eea
 From Proposition (\ref{pro1}), the $(C,C')-$controlled g-frame
 operator $S_{(C,C')}$ is positive, self adjoint and invertible.
 Hence
 \bea\label{eqn7}
 \big\langle S_{(c,c')}^{\frac{1}{2}}f,
 S_{(c,c')}^{\frac{1}{2}}f\big\rangle=\big\langle S_{(c,c')}f,
 f\big\rangle=\sum_{j\in J}\langle
\Lambda_{j}Cf,\Lambda_{j}C'f\rangle.
 \eea
 Using (\ref{eqn7}) in (\ref{eqn6}), we get
 \bea\label{eqn8}
 \sqrt{A}\|f\|\leq \|
 S_{(c,c')}^{\frac{1}{2}}f\|\leq\sqrt{B}\|f\|,~~\forall f\in \mathcal{H}.
 \eea
 According to Lemma \ref{l3.1} and inequality (\ref{eqn8}), there exist
 constant $m,M> 0$ such that
 \bee
 m\langle f,f\rangle\leq \big\langle S_{(c,c')}^{\frac{1}{2}}f,
 S_{(c,c')}^{\frac{1}{2}}f\big\rangle=\sum_{j\in J}\langle
\Lambda_{j}Cf,\Lambda_{j}C'f\rangle\leq M\langle
f,f\rangle,~~\forall f\in \mathcal{H}.
 \eee
 Therefore, $\{\Lambda_{j}:j\in J\}$ is a $(C,C')-$controlled
 g-frame for $\mathcal{H}$ with respect to $\{\mathcal{H}_{j}\}_{j\in J}$.
\end{proof}

\begin{d1}\rm
Let $C\in GL^{+}(\mathcal{H})$. The sequence $\{\Lambda_{j}\in
End_{A}^{*}(\mathcal{H},\mathcal{H}_{j}):j\in J\}$ is said to be a
$(C,C)-$controlled g-frame or $C^{2}-$controlled g-frame if there
exist constants $0<A\leq B< \infty$ such that \bee A\langle
f,f\rangle\leq\sum_{j\in J}\langle \Lambda_{j}C f, \Lambda_{j}C
 f\rangle\leq B\langle f,f\rangle,~\forall f\in \mathcal{H},\eee
 or equivalently,
 \bea\label{eqn9} A
\|f\|^{2}\leq\|\sum_{j\in J}\langle \Lambda_{j}C f, \Lambda_{j}C
 f\rangle\|\leq B\|f\|^{2},~\forall f\in \mathcal{H}.\eea
\end{d1}

 Using some
tools from operator algebras, Xiao and Zeng \cite{xiao} have
proved the following equivalent characterization of g-frames in
Hilbert $C^{*}$-modules.
\begin{t1}\rm \label{t2}\cite{xiao}
Let $\{ \Lambda_j \in
End_{A}^{*}(\mathcal{H},{\mathcal{H}}_j):j\in J\}$ and $\sum_{j\in
J}\langle \Lambda_{j}f,\Lambda_{j}f\rangle$ converge in norm for
$f\in \mathcal{H}$. Then $\{ \Lambda_j :j\in J\}$ is a g-frame for
$\mathcal{H}$ with respect to $\{\mathcal{H}_{j}\}_{j\in J}$ if
and only if there exist constants $A,B> 0$ such that \bea
\label{eqn16} A\|f\|^{2}\leq \big\|\sum_{j\in J}\langle
\Lambda_{j}f,\Lambda_{j}f\rangle\big\|\leq B\|f\|^{2},~~\forall
f\in \mathcal{H}. \eea
\end{t1}
\noindent The above result can be easily seen as a corollary of
our Theorem \ref{t1},  when we take $C=C'=I$.

\begin{p1}\rm
Let $C\in GL^{+}(\mathcal{H})$. The family $\{\Lambda_{j}:j\in
J\}$ is a g-frame if and only if $\{\Lambda_{j}:j\in J\}$ is a
$C^{2}-$controlled g-frame.
\end{p1}

\begin{proof}
Suppose that $\{\Lambda_{j}:j\in J\}$ is a $C^{2}-$controlled
g-frame with bounds $A$ and $B$. Then from (\ref{eqn9}), we have
\bee A \|f\|^{2}\leq\|\sum_{j\in J}\langle \Lambda_{j}C f,
\Lambda_{j}C
 f\rangle\|\leq B\|f\|^{2},~\forall f\in \mathcal{H}.\eee
 Now for any $f\in \mathcal{H}$,\bee
 A\|f\|^{2}=A\|CC^{-1}f\|^{2}&\leq& A\|C\|^{2}\|C^{-1}f\|^{2}\\
 &\leq& \|C\|^{2}\|\sum_{j\in J}\langle \Lambda_{j}CC^{-1}f,
 \Lambda_{j}CC^{-1}f\rangle\|\\
 &=& \|C\|^{2}\|\sum_{j\in J}\langle \Lambda_{j}f,
 \Lambda_{j}f\rangle\|.
 \eee
 Hence
 \bea\label{eqn10}
 A\|C\|^{-2}\|f\|^{2}\leq \|\sum_{j\in J}\langle \Lambda_{j}f,
 \Lambda_{j}f\rangle\|.
 \eea
 Again for any $f\in \mathcal{H}$,
 \bea\label{eqn11}
\|\sum_{j\in J}\langle \Lambda_{j}f,
 \Lambda_{j}f\rangle\|&=&\|\sum_{j\in J}\langle \Lambda_{j}CC^{-1}f,
 \Lambda_{j}CC^{-1}f\rangle\|\nn\\
 &\leq& B\|C^{-1}f\|^{2}\leq B\|C^{-1}\|^{2}\|f\|^{2}.
 \eea
 From (\ref{eqn10}), (\ref{eqn11}) and Theorem \ref{t2}, we conclude that $\{\Lambda_{j}:j\in J\}$ is a
g-frame with bound $A\|C\|^{-2}$ and $B\|C^{-1}\|^{2}$.\\
Conversely, let $\{\Lambda_{j}:j\in J\}$ is a g-frame with bounds
$A'$ and $B'$. Then for all $f\in \mathcal{H}$, \bee A'\langle
f,f\rangle\leq\sum_{j\in J}\langle \Lambda_{j} f, \Lambda_{j}
 f\rangle\leq B'\langle f,f\rangle.
\eee So for $f\in \mathcal{H}$ we have $Cf\in\mathcal{H}$, and
 \bea\label{eqn12} \sum_{j\in J}\langle \Lambda_{j}C f,
\Lambda_{j}C
 f\rangle\leq B'\langle C f,Cf\rangle\leq B'\|C\|^{2}\langle
 f,f\rangle.
\eea Also for any $f\in \mathcal{H}$, \bea\label{eqn13} A'\langle
f,f\rangle= A'\langle C^{-1}C f,C^{-1}C f\rangle&\leq& A'
\|C^{-1}\|^{2}\langle Cf,Cf\rangle\nn\\
&\leq& \|C^{-1}\|^{2}\sum_{j\in J}\langle \Lambda_{j}C
f,\Lambda_{j}Cf\rangle.\eea From (\ref{eqn12}) and (\ref{eqn13}),
we have \bee A' \|C^{-1}\|^{-2}\langle f,f\rangle\leq  \sum_{j\in
J}\langle \Lambda_{j}C f,\Lambda_{j}Cf\rangle\leq
B'\|C\|^{2}\langle
 f,f\rangle,~~\forall f\in \mathcal{H}.\eee
 Hence $\{\Lambda_{j}:j\in J\}$ is a
$C^{2}-$controlled g-frame with bounds $A' \|C^{-1}\|^{-2}$ and
$B'\|C\|^{2}$.
\end{proof}
Next, we study when a $g$-frame becomes a $(C,C')-$controlled
$g$-frame.
\begin{p1}\rm
Assume that $\{\Lambda_{j}:j\in J\}$ is a g-frame for the Hilbert
$C^{*}$-module $\mathcal{H}$ with respect to
$\{\mathcal{H}_{j}\}_{j\in J}$. Let $S_{\Lambda}$ be the g-frame
operator associated with the g-frame $\{\Lambda_{j}:j\in J\}$ as
defined in (\ref{eqn18}). Let $C,C'\in GL^{+}(\mathcal{H})$. Then
$\{\Lambda_{j}:j\in J\}$ is a $(C,C')-$controlled g-frame.
\end{p1}
\begin{proof}
$\{\Lambda_{j}:j\in J\}$ is a g-frame for the Hilbert
$C^{*}$-module $\mathcal{H}$ with bounds $A$ and $B$. By the
equivalence condition (\ref{eqn16}) of g-frame, we have
\bea\label{eqn14} &&A\|f\|^{2} \leq \big\| \sum_{j\in J}\langle
\Lambda_{j}f, \Lambda_{j}f\rangle\big\|\leq B\|f\|^{2},
~~\forall f\in \mathcal{H}\nn\\
&\Rightarrow& A\|f\|^{2}\leq \big\|\langle S_{\Lambda}f,f\rangle
\big\|\leq B\|f\|^{2},~~\forall f\in \mathcal{H}. \eea Again we
have \bee \big\| \sum_{j\in J}\langle \Lambda_{j}C f, \Lambda_{j}
C'f\rangle\big\|=\big\|\langle S_{C,C'}f,f\rangle \big\| \eee and
\bea\label{eqn15} \big\| \sum_{j\in J}\langle \Lambda_{j}C f,
\Lambda_{j} C'f\rangle\big\| =\|C\|\|C'\|\big\| \sum_{j\in
J}\langle \Lambda_{j}f,
\Lambda_{j}f\rangle\big\|=\|C\|\|C'\|\big\|\langle
S_{\Lambda}f,f\rangle\big\|.\eea From (\ref{eqn14}) and
(\ref{eqn15}), we have \bee A\|C\|\|C'\| \|f\|^{2}\leq \big\|
\sum_{j\in J}\langle \Lambda_{j}C f, \Lambda_{j}
C'f\rangle\big\|\leq B\|C\|\|C'\|\|f\|^{2},~~\forall f\in
\mathcal{H}.\eee By using Theorem \ref{t1}, we conclude that
$\{\Lambda_{j}:j\in J\}$ is a $(C,C')-$controlled g-frame with
bounds $A\|C\|\|C'\|$ and $B\|C\|\|C'\|$.
\end{proof}

\begin{t1}\label{t3}\rm
Let $C,C'\in GL^{+}(\mathcal{H})$, $\{\Lambda_{j}:j\in J\}\subset
End_{A}^{*}(\mathcal{H},\mathcal{H}_{j})$, and $C, C'$ commute
with each other and commute with $\Lambda_{j}^{*}\Lambda_{j}$
 for all $j\in J$. Then the sequence
$\{\Lambda_{j}:j\in J\}$ is a $(C,C')-$controlled $g$-Bessel
sequence for $\mathcal{H}$ with respect to
$\{\mathcal{H}_{j}\}_{j\in J}$ with bound $B$ if and only if the
operator $T_{(C,C^{\prime})}:l^{2}(\{\mathcal{H}_{j}\}_{j\in
J})\rightarrow \mathcal{H}$ given by \bee
T_{(C,C^{\prime})}(\{g_{j}\}_{j\in J})=\sum_{j\in J}
(CC')^{\frac{1}{2}}\Lambda_{j}^{*}g_{j},~~\forall \{g_{j}\}_{j\in
J}\in l^{2}(\{\mathcal{H}_{j}\}_{j\in J})\eee is well defined and
bounded operator with $\|T_{(C,C^{\prime})}\|\leq \sqrt{B}$.
\end{t1}
\begin{proof}
Let $\{\Lambda_{j}:j\in J\}$ be a $C,C'-$controlled g-Bessel
sequence for $\mathcal{H}$ with respect to
$\{\mathcal{H}_{j}\}_{j\in J}$ with bound $B$. As a result of
Theorem \ref{t1}, \bea\label{eqn23}\|\sum_{j\in J}\langle
\Lambda_{j}Cf,\Lambda_{j}C'f\rangle\|\leq B\|f\|^{2},~~\forall
f\in \mathcal{H}.
 \eea

 For any sequence $\{g_{j}\}_{j\in J}\in l^{2}(\{\mathcal{H}_{j}\}_{j\in
 J})$,
 \bee
 \big\|T_{(C,C')}\big(\{g_{j}\}_{j\in J}\big)\big\|^{2}&=&\sup_{f\in \mathcal{H},~\|f\|=1}
 \big\|\langle T_{(C,C')}\big(\{g_{j}\}_{j\in J}\big),f\rangle\big\|^{2}\\
 &=&\sup_{f\in \mathcal{H},~\|f\|
 =1}\big\|\big\langle \sum_{j\in
 J}(CC')^\frac{1}{2}\Lambda_{j}^{*}g_{j},f\big\rangle\big\|^{2}\\
 &=& \sup_{f\in \mathcal{H},~\|f\|
 =1}\big\| \sum_{j\in
 J}\big\langle(CC')^\frac{1}{2}\Lambda_{j}^{*}g_{j},f\big\rangle\big\|^{2}\\
 &=& \sup_{f\in \mathcal{H},~\|f\|
 =1}\big\| \sum_{j\in
 J}\big\langle g_{j},
 \Lambda_{j}(CC')^\frac{1}{2}f\big\rangle\big\|^{2}\\
 &\leq& \sup_{f\in \mathcal{H},~\|f\|
 =1}\big\|\sum_{j\in J}\langle g_{j},g_{j}\rangle\big\|\big\| \sum_{j\in
 J}\big\langle  \Lambda_{j}(CC')^\frac{1}{2}f,
 \Lambda_{j}(CC')^\frac{1}{2}f\big\rangle\big\|\\
 &=& \sup_{f\in \mathcal{H},~\|f\|
 =1}\big\|\sum_{j\in J}\langle g_{j},g_{j}\rangle\big\|\big\| \sum_{j\in
 J}\big\langle  \Lambda_{j}C f,
 \Lambda_{j}C'f\big\rangle\big\|\\
& \leq& \sup_{f\in \mathcal{H},~\|f\|
 =1}\big\|\sum_{j\in J}\langle g_{j},g_{j}\rangle\big\|
 B\|f\|^{2}=B \big\|\{g_{j}\}\big\|^{2}.
 \eee
 Therefore, the sum $\displaystyle\sum_{j\in J}
(CC')^{\frac{1}{2}}\Lambda_{j}^{*}g_{j}$ is convergent, and we
have \bee &&\big\|T_{(C,C^{\prime})}(\{g_{j}\}_{j\in
J})\big\|^{2}\leq B\big\|\{g_{j}\}\big\|^{2}\\
&\Rightarrow& \|T_{(C,C')}\|\leq \sqrt{B}.\eee Hence the operator
$T_{(C,C')}$ is well defined, bounded and $\|T_{(C,C')}\|\leq
\sqrt{B}$.\\

Conversely, let the operator $T_{(C,C')}$ is well defined, bounded
and $\|T_{(C,C')}\|\leq \sqrt{B}$.\\ For any $f\in \mathcal{H}$
and finite subset $K\subset J$, we have \bee \sum_{j\in
K}\big\langle \Lambda_{j}C f,\Lambda_{j}C' f \big\rangle &=&
\sum_{j\in
K}\big\langle C'\Lambda_{j}^{*}\Lambda_{j}C f, f \big\rangle\\
&=&  \sum_{j\in K}\big\langle
(CC')^{\frac{1}{2}}\Lambda_{j}^{*}\Lambda_{j}(CC')^{\frac{1}{2}}
f, f \big\rangle\\
&=& \Big\langle T_{(C,C^{\prime})}(\{g_{j}\}_{j\in
J}),f\Big\rangle\\
&\leq& \|T_{(C,C^{\prime})}\|\|\{g_{j}\}_{j\in J}\|\|f\|,\eee
where \bee
    g_{j}=
    \begin{cases}
      \Lambda_{j}(CC')^{\frac{1}{2}}f, & \text{if}\ j\in K \\
      0, & j \notin K
    \end{cases}.
  \eee
  Therefore,
  \bee
\sum_{j\in K}\big\langle \Lambda_{j}C f,\Lambda_{j}C' f
\big\rangle&\leq& \|T_{(C,C^{\prime})}\|\Big(\sum_{j\in
K}\|\Lambda_{j}(CC')^{\frac{1}{2}}
f\|^{2}\Big)^{\frac{1}{2}}\|f\|\\
&=& \|T_{(C,C^{\prime})}\|\Big(\sum_{j\in K}\big\langle
\Lambda_{j}C f, \Lambda_{j}C' f
\big\rangle\Big)^{\frac{1}{2}}\|f\|.
  \eee
Since $K$ is arbitrary, we have
  \bee
&& \sum_{j\in J}\big\langle \Lambda_{j}C f,\Lambda_{j}C' f
\big\rangle \leq \|T_{(C,C^{\prime})}\|^{2}\|f\|^{2}\\
&\Rightarrow& \sum_{j\in J}\big\langle \Lambda_{j}C
f,\Lambda_{j}C' f\big\rangle \leq B\|f\|^{2}, {~\rm as~}
\|T_{(C,C^{\prime})}\|\leq \sqrt{B}.
  \eee
Hence we conclude that $\{\Lambda_{j}:j\in J\}$ is a
$(C,C')-$controlled g-Bessel sequence for $\mathcal{H}$ with
respect to $\{\mathcal{H}_{j}\}_{j\in J}$.
\end{proof}
Now we prove some perturbation results for $(C,C')-$controlled
g-frame.
\begin{t1}\rm
Let $\{\Lambda_{j}\in End_{A}^{*}(\mathcal{H},\mathcal{H}_j):j\in
J\}$ be a $(C,C')-$controlled g-frame for $\mathcal{H}$ with
respect to $\{\mathcal{H}_j\}_{j\in J}$. Let $\{\Pi_{j}\in
End_{A}^{*}(\mathcal{H},\mathcal{H}_j):j\in J\}$ be any sequence,
and assume that $C$ and $C'$ commute with each other and commute
with $(\Lambda_{j}-\Pi_{j})^{*}(\Lambda_{j}-\Pi_{j})$. Then
$\{\Pi_{j}:j\in J\}$ is a $(C,C')-$controlled g-frame for
$\mathcal{H}$ with respect to $\{\mathcal{H}_j\}_{j\in J}$ if and
only if there exists constants $M_{1}$ and $M_{2}$ such that
\bea\label{eqn19} \big\|\sum_{j\in J}\big\langle
(\Lambda_{j}-\Pi_{j})C f, (\Lambda_{j}-\Pi_{j})C' f
\big\rangle\big\|\leq M_{1}\big\|\sum_{j\in J}\big\langle
\Lambda_{j}C f, \Lambda_{j}C' f \big\rangle\big\|\\
{\rm and~}\label{eqn20} \big\|\sum_{j\in J}\big\langle
(\Lambda_{j}-\Pi_{j})C f, (\Lambda_{j}-\Pi_{j})C' f
\big\rangle\big\|\leq M_{2}\big\|\sum_{j\in J}\big\langle \Pi_{j}C
f, \Pi_{j}C' f \big\rangle\big\|.\eea
\end{t1}
\begin{proof}
Let $\{\Lambda_{j}:j\in J\}$ be a $(C,C')-$controlled g-frame for
$\mathcal{H}$ with lower and upper bounds $A_{1}$ and $B_{1}$,
respectively. Also suppose that $\{\Pi_{j}:j\in J\}$ be a
$(C,C')-$controlled g-frame for $\mathcal{H}$ with lower and upper
bounds $A_{2}$ and $B_{2}$, respectively. Then \bee
&&\big\|\sum_{j\in J}\big\langle (\Lambda_{j}-\Pi_{j})C f,
(\Lambda_{j}-\Pi_{j})C' f \big\rangle\big\|\\&=&\big\|\sum_{j\in
J}\big\langle (\Lambda_{j}-\Pi_{j})(CC')^{\frac{1}{2}} f,
(\Lambda_{j}-\Pi_{j})(CC')^{\frac{1}{2}} f \big\rangle\big\|
\\
&=& \big\|\big\{(\Lambda_{j}-\Pi_{j})(CC')^{\frac{1}{2}}
f\big\}_{j\in J}\big\|^{2}\\
&\leq& \big\|\big\{\Lambda_{j}(CC')^{\frac{1}{2}} f\big\}_{j\in
J}\big\|^{2}+\big\|\big\{\Pi_{j}(CC')^{\frac{1}{2}} f\big\}_{j\in
J}\big\|^{2}\\
&=& \big\|\sum_{j\in J}\big\langle \Lambda_{j}(CC')^{\frac{1}{2}}
f, \Lambda_{j}(CC')^{\frac{1}{2}} f
\big\rangle\big\|+\big\|\sum_{j\in J}\big\langle
\Pi_{j}(CC')^{\frac{1}{2}} f, \Pi_{j}(CC')^{\frac{1}{2}} f
\big\rangle\big\|\\
&=& \big\|\sum_{j\in J}\big\langle \Lambda_{j}C f, \Lambda_{j}C' f
\big\rangle\big\|+\big\|\sum_{j\in J}\big\langle \Pi_{j}C f,
\Pi_{j}C' f \big\rangle\big\|\\
&\leq&  \big\|\sum_{j\in J}\big\langle \Lambda_{j}C f,
\Lambda_{j}C' f \big\rangle\big\|+B_{2}\|f\|^{2}\\
&\leq& \big\|\sum_{j\in J}\big\langle \Lambda_{j}C f,
\Lambda_{j}C' f
\big\rangle\big\|+\frac{B_{2}}{A_{1}}\big\|\sum_{j\in
J}\big\langle \Lambda_{j}C f, \Lambda_{j}C' f
\big\rangle\big\|\\
&=& \Big(1+\frac{B_{2}}{A_{1}}\Big)\big\|\sum_{j\in J}\big\langle
\Lambda_{j}C f, \Lambda_{j}C' f \big\rangle\big\|.\eee Thus
(\ref{eqn19}) is proved, where
$M_{1}=\Big(1+\frac{B_{2}}{A_{1}}\Big)$. In a similar manner, one
can obtain \bee \big\|\sum_{j\in J}\big\langle
(\Lambda_{j}-\Pi_{j})C f, (\Lambda_{j}-\Pi_{j})C' f
\big\rangle\big\|\leq
\Big(1+\frac{B_{1}}{A_{2}}\Big)\big\|\sum_{j\in J}\big\langle
\Pi_{j}C f, \Pi_{j}C' f \big\rangle\big\|.\eee Hence (\ref{eqn20})
follows with $M_{2}=\Big(1+\frac{B_{1}}{A_{2}}\Big)$.\\

Conversely, suppose that $\{\Lambda_{j}:j\in J\}$ be a
$(C,C')-$controlled g-frame for $\mathcal{H}$ with lower and upper
bounds $A_{1}$ and $B_{1}$, respectively, and (\ref{eqn19}) and
(\ref{eqn20}) hold true. Then for any $f\in \mathcal{H}$, using
(\ref{eqn20}) we get \bea A_{1}\|f\|^{2}&\leq&  \big\|\sum_{j\in
J}\big\langle \Lambda_{j}C f,
\Lambda_{j}C' f \big\rangle\big\|\nn\\
&=&\big\|\sum_{j\in J}\big\langle \Lambda_{j}(CC')^{\frac{1}{2}}
f, \Lambda_{j}(CC')^{\frac{1}{2}} f \big\rangle\big\|\nn\\
&=&  \big\|\big\{\Lambda_{j}(CC')^{\frac{1}{2}} f\big\}_{j\in
J}\big\|^{2}\nn\\
&\leq& \big\|\big\{(\Lambda_{j}-\Pi_{j})(CC')^{\frac{1}{2}}
f\big\}_{j\in J}\big\|^{2}+\big\|\big\{\Pi_{j}(CC')^{\frac{1}{2}}
f\big\}_{j\in J}\big\|^{2}\nn\\
&=&\big\|\sum_{j\in J}\big\langle
(\Lambda_{j}-\Pi_{j})(CC')^{\frac{1}{2}} f,
(\Lambda_{j}-\Pi_{j})(CC')^{\frac{1}{2}} f
\big\rangle\big\|\nn\\&&+\big\|\sum_{j\in J}\big\langle
\Pi_{j}(CC')^{\frac{1}{2}} f, \Pi_{j}(CC')^{\frac{1}{2}} f
\big\rangle\big\|\nn\\
&=& \big\|\sum_{j\in J}\big\langle (\Lambda_{j}-\Pi_{j})C f,
(\Lambda_{j}-\Pi_{j})C' f \big\rangle\big\|+\big\|\sum_{j\in
J}\big\langle \Pi_{j}C f, \Pi_{j}C' f \big\rangle\big\|\nn\\
&\leq& M_{2}\big\|\sum_{j\in J}\big\langle \Pi_{j}C f, \Pi_{j}C' f
\big\rangle\big\|+\big\|\sum_{j\in J}\big\langle \Pi_{j}C f,
\Pi_{j}C' f \big\rangle\big\|\nn\\&=&(1+M_{2})\big\|\sum_{j\in
J}\big\langle \Pi_{j}C f, \Pi_{j}C' f \big\rangle\big\|.\nn \eea
This implies that \bea\label{eqn21} \frac{A_{1}}{
(1+M_{2})}\|f\|^{2}\leq\big\|\sum_{j\in J}\big\langle \Pi_{j}C f,
\Pi_{j}C' f \big\rangle\big\|.\eea Also we have \bea\label{eqn22}
\big\|\sum_{j\in J}\big\langle \Pi_{j}C f, \Pi_{j}C' f
\big\rangle\big\|&=&\big\|\sum_{j\in J}\big\langle
\Pi_{j}(CC')^{\frac{1}{2}} f, \Pi_{j}(CC')^{\frac{1}{2}} f
\big\rangle\big\|\nn\\
&=& \big\|\big\{\Pi_{j}(CC')^{\frac{1}{2}} f\big\}_{j\in
J}\big\|^{2}\nn\\
&\leq& \big\|\big\{\Lambda_{j}(CC')^{\frac{1}{2}} f\big\}_{j\in
J}\big\|^{2}+\big\|\big\{(\Lambda_{j}-\Pi_{j})(CC')^{\frac{1}{2}}
f\big\}_{j\in J}\big\|^{2}\nn\\
&=& \big\|\sum_{j\in J}\big\langle
(\Lambda_{j}-\Pi_{j})(CC')^{\frac{1}{2}} f,
(\Lambda_{j}-\Pi_{j})(CC')^{\frac{1}{2}} f
\big\rangle\big\|\nn\\&&+\big\|\sum_{j\in J}\big\langle
\Lambda_{j}(CC')^{\frac{1}{2}} f, \Lambda_{j}(CC')^{\frac{1}{2}} f
\big\rangle\big\|\nn\\
&=& \big\|\sum_{j\in J}\big\langle (\Lambda_{j}-\Pi_{j})C f,
(\Lambda_{j}-\Pi_{j})C' f \big\rangle\big\|+\big\|\sum_{j\in
J}\big\langle \Lambda_{j}C f, \Lambda_{j}C' f
\big\rangle\big\|\nn\\
&\leq& M_{1}\big\|\sum_{j\in J}\big\langle \Lambda_{j}C f,
\Lambda_{j}C' f \big\rangle\big\|+\big\|\sum_{j\in J}\big\langle
\Lambda_{j}C f, \Lambda_{j}C' f \big\rangle\big\|\nn\\
&=&(1+M_{1})\big\|\sum_{j\in J}\big\langle \Lambda_{j}C f,
\Lambda_{j}C' f \big\rangle\big\|\nn\\
&\leq& (1+M_{1})B_{1}\|f\|^{2}.\eea Therefore from (\ref{eqn21})
and (\ref{eqn22}), it is clear that $\{\Pi_{j}:j\in J\}$ is a
$(C,C')-$controlled g-frame for $\mathcal{H}$ with respect to
$\{\mathcal{H}_{j}\}_{j\in J}$.
\end{proof}
\begin{p1}\rm
Let $\{\Lambda_{j}\in End_{A}^{*}(\mathcal{H},\mathcal{H}_j):j\in
J\}$ and $\{\Gamma_{j}\in
End_{A}^{*}(\mathcal{H},\mathcal{H}_j):j\in J\}$ be two
$(C,C')$-controlled $g$-Bessel sequences for $\mathcal{H}$ with
respect to $\{\mathcal{H}\}_{j\in J}$ with bounds $B_{1}$ and
$B_{2}$, respectively. Then the operator
$L_{(C,C')}:\mathcal{H}\rightarrow \mathcal{H}$ given by
\bea\label{eqn24} L_{(C,C')}(f)=\displaystyle\sum_{j\in
J}C'\Gamma_{j}^{*}\Lambda_{j}Cf\eea is well defined and bounded
with $\|L_{(C,C')}\|\leq \sqrt{B_{1}B_{2}}$. Also its adjoint
operator is $L_{(C,C')}^{*}(g)=\displaystyle\sum_{j\in
J}C\Lambda_{j}^{*}\Gamma_{j}C'g$.
\end{p1}
\begin{proof}
For any $f\in\mathcal{H}$ and $K\subset J$, we have

\bee \big\|\sum_{j\in
K}C'\Gamma_{j}^{*}\Lambda_{j}Cf\big\|^{2}&=&\sup_{g\in
\mathcal{H},\|g\|=1}\big\|\Big\langle \sum_{j\in
K}C'\Gamma_{j}^{*}\Lambda_{j}Cf,~g \Big\rangle\big\|^{2}\\
&=& \sup_{g\in \mathcal{H},\|g\|=1}\big\|\sum_{j\in K}\Big\langle
\Lambda_{j}Cf,\Gamma_{j}C'g \Big\rangle\big\|^{2}\\
&\leq& \sup_{g\in \mathcal{H},\|g\|=1}\big\|\sum_{j\in
K}\Big\langle \Lambda_{j}Cf,\Lambda_{j}Cf \Big\rangle\big\|
\big\|\sum_{j\in K}\Big\langle \Gamma_{j}C'g,\Gamma_{j}C'g
\Big\rangle\big\|\\
&\leq& \big\|\sum_{j\in K}\Big\langle \Lambda_{j}Cf,\Lambda_{j}Cf
\Big\rangle\big\|~ B_{2}\\
&\leq& B_{1}B_{2}\|f\|^{2}.\eee Since $K$ is arbitrary the series
$\displaystyle\sum_{j\in J}C'\Gamma_{j}^{*}\Lambda_{j}Cf$
converges in $\mathcal{H}$, and \bee
\|L_{(C,C')}\|=\big\|\sum_{j\in
K}C'\Gamma_{j}^{*}\Lambda_{j}Cf\big\|\leq \sqrt{B_{1}B_{2}}.\eee
Moreover, we see that \bee\langle L_{(C,C')}f,g\rangle=\Big\langle
\sum_{j\in K}C'\Gamma_{j}^{*}\Lambda_{j}Cf,~g
\Big\rangle=\sum_{j\in K}\Big\langle
C'\Gamma_{j}^{*}\Lambda_{j}Cf,~g \Big\rangle &=&\sum_{j\in
K}\Big\langle f~,C\Lambda_{j}^{*}\Gamma_{j}C'g \Big\rangle\\
&=& \Big\langle f~,\sum_{j\in K}C\Lambda_{j}^{*}\Gamma_{j}C'g
\Big\rangle.\eee
 Thus $L_{(C,C')}^{*}(g)=\displaystyle\sum_{j\in
J}C\Lambda_{j}^{*}\Gamma_{j}C'g$.
\end{proof}
\begin{t1}\rm
Let $\{\Lambda_{j}\in End_{A}^{*}(\mathcal{H},\mathcal{H}_j):j\in
J\}$ be a $(C,C')$-controlled $g$-frame for $\mathcal{H}$ with
respect to $\{\mathcal{H}\}_{j\in J}$, and  $\{\Gamma_{j}\in
End_{A}^{*}(\mathcal{H},\mathcal{H}_j):j\in J\}$ be a
$(C,C')$-controlled $g$-Bessel sequence for $\mathcal{H}$ with
respect to $\{\mathcal{H}\}_{j\in J}$. Assume that $C$ and $C'$
commute with each other and commute with
$\Gamma_{j}^{*}\Gamma_{j}$. If the operator $L_{(C,C')}$ defined
in (\ref{eqn24}) is surjective then $\{\Gamma_{j}:j\in J\}$ is
also a $(C,C')$-controlled $g$-frame for $\mathcal{H}$ with
respect to $\{\mathcal{H}\}_{j\in J}$.
\end{t1}

\begin{proof}
It is given that $\{\Lambda_{j}:j\in J\}$ is a $(C,C')$-controlled
$g$-frame for $\mathcal{H}$ with respect to $\{\mathcal{H}\}_{j\in
J}$. Then by Theorem \ref{t3}, the operator
$T_{(C,C^{\prime})}:l^{2}(\{\mathcal{H}_{j}\}_{j\in J})\rightarrow
\mathcal{H}$ given by \bee T_{(C,C^{\prime})}(\{g_{j}\}_{j\in
J})=\sum_{j\in J}
(CC')^{\frac{1}{2}}\Lambda_{j}^{*}g_{j},~~\forall \{g_{j}\}_{j\in
J}\in l^{2}(\{\mathcal{H}_{j}\}_{j\in J})\eee is well defined and
bounded operator. By (\ref{eqn3}) its adjoint operator
$T^{*}_{(C,C^{\prime})}: \mathcal{H}\rightarrow
l^{2}(\{\mathcal{H}_{j}\}_{j\in J})$ is given by \bee
T^{*}_{(C,C^{\prime})}(f)=\big\{\Lambda_{j}(C^{\prime}C)^{\frac{1}{2}}f\big\}_{j\in
J},~\forall f\in \mathcal{H}. \eee Since $\{\Gamma_{j}:j\in J\}$
is also a $(C,C')$-controlled $g$-Bessel sequence for
$\mathcal{H}$ with respect to $\{\mathcal{H}\}_{j\in J}$, again by
Theorem \ref{t3}, the operator
$P_{(C,C^{\prime})}:l^{2}(\{\mathcal{H}_{j}\}_{j\in J})\rightarrow
\mathcal{H}$ given by \bee P_{(C,C^{\prime})}(\{g_{j}\}_{j\in
J})=\sum_{j\in J} (CC')^{\frac{1}{2}}\Gamma_{j}^{*}g_{j},~~\forall
\{g_{j}\}_{j\in J}\in l^{2}(\{\mathcal{H}_{j}\}_{j\in J})\eee is
well defined and bounded operator. Again its adjoint operator is
given by \bee
P^{*}_{(C,C^{\prime})}(f)=\big\{\Gamma_{j}(C^{\prime}C)^{\frac{1}{2}}f\big\}_{j\in
J},~\forall f\in \mathcal{H}.\eee Hence for any $f\in
\mathcal{H}$, the operator defined in (\ref{eqn24}) can be written
as \bee L_{(C,C')}(f)=\displaystyle\sum_{j\in
J}C'\Gamma_{j}^{*}\Lambda_{j}Cf=P_{(C,C')}T^{*}_{(C,C')}(f).\eee
Since $L_{(C,C')}$ is surjective then for any $f\in \mathcal{H}$,
there exists $g\in \mathcal{H}$ such that
$f=L_{(C,C')}(g)=P_{(C,C')}T^{*}_{(C,C')}(g)$, and
$T^{*}_{(C,C')}(g)\in l^{2}(\{\mathcal{H}_{j}\}_{j\in J})$. This
implies that $P_{(C,C')}$ is surjective. As a result of Lemma
\ref{l3.1}, we have $P^{*}_{(C,C')}$ is bounded below, that is
there exists $m>0$ such that \bee &&\big\langle
P^{*}_{(C,C')}f,P^{*}_{(C,C')}f\big\rangle\geq m\langle
f,f\rangle,\forall f\in \mathcal{H}\nn\\
&\Rightarrow& \big\langle
P_{(C,C')}P^{*}_{(C,C')}f,f\big\rangle\geq m\langle
f,f\rangle,\forall f\in \mathcal{H}\nn \\
&\Rightarrow& \Big\langle \sum_{j\in
J}(CC')^{\frac{1}{2}}\Gamma_{j}^{*}\Gamma_{j}(C'C)^{\frac{1}{2}}f,~f\Big\rangle\geq
m\langle f,f\rangle,~\forall
f\in \mathcal{H}\nn\\
&\Rightarrow& \sum_{j\in J}\Big\langle \Gamma_{j}C f,\Gamma_{j}C'
f\Big\rangle\geq m\langle f,f\rangle,~\forall f\in \mathcal{H}.
\eee Hence $\{\Gamma_{j}:j\in J\}$ is also a $(C,C')$-controlled
$g$-frame for $\mathcal{H}$ with respect to $\{\mathcal{H}\}_{j\in
J}$.
\end{proof}
 
\end{document}